\input amstex\documentstyle{amsppt}  
\pagewidth{12.5cm}\pageheight{19cm}\magnification\magstep1  
\topmatter
\title Adjacency for special representations of a Weyl group\endtitle
\author G. Lusztig\endauthor
\address{Department of Mathematics, M.I.T., Cambridge, MA 02139}\endaddress
\thanks{Supported by NSF grant DMS-1855773.}\endthanks
\endtopmatter   
\document
\define\tcu{\ti{\cu}}
\define\sgn{\text{\rm sgn}}
\define\hcs{\hat{\cs}}

\define\Irr{\text{\rm Irr}}

\define\si{\sim}

\define\qua{\quad}

\define\bC{\bar C}

\define\lb{\linebreak}

\define\part{\partial}
\define\emp{\emptyset}

\define\n{\notin}

\define\m{\mapsto}
\define\do{\dots}

\define\esm{\endsmallmatrix}
\define\sub{\subset}    

\define\T{\times}
\define\ti{\tilde}
\define\nl{\newline}

\define\ot{\otimes}

\define\a{\alpha}

\define\e{\epsilon}

\define\io{\iota}

\define\x{\xi}

\redefine\D{\Delta}

\define\CC{\bold C}

\define\NN{\bold N}

\define\QQ{\bold Q}

\define\ZZ{\bold Z}

\define\cl{\Cal L}

\define\cs{\Cal S}

\define\cu{\Cal U}

\define\fS{\frak S}
\define\fT{\frak T}

\define\bul{\bullet}

\define\cir{\circ}

\head Introduction\endhead
\subhead 0.1\endsubhead
Let $G$ be a connected reductive group over $\CC$. Let $\cu_G$ be the set
of unipotent conjugacy classes of $G$. Let $C'\in\cu_G$ and let $C\in\cu_G$
be maximal with the property that $C\sub\bC'-C'$ ($\bC'$ is the closure of $C'$).
A remarkable result of Kraft and Procesi \cite{KP81}, \cite{KP82} is that
when $G$ is a classical group, then either $\dim(C')=\dim(C)+2$ or the singularity
of $C'$ at a point of $C$ is the same as the singularity at $1$ of a minimal unipotent
class in a smaller reductive group. This result has been recently extended
to exceptional groups by Fu, Juteau, Levy and Sommers \cite{FJLS}.

In the late 1980's, inspired by \cite{KP81}, \cite{KP82}, I showed (unpublished) that the
results of {\it loc.cit.} have a (weak) analogue in the case where $\cu_G$ is replaced
by $\cu_G^{sp}$ (the set of special unipotent classes of $G$). The analogues in this
case of the pairs $C,C'$ above can be viewed as edges of a graph with set of vertices
$\cu_G^{sp}$. One feature that was not present in {\it loc.cit.} (except in type $A$)
is that $\cu_G^{sp}$ has an order reversing involution which preserves the graph structure
and that if two edges of the graph are interchanged by this involution then
at least one of them is associated to a pair $C,C'$ with $\dim(C')=\dim(C)\pm2$.
(see Theorem 0.3).

This property allows us to construct the graph above purely in terms of
the involution above and truncated induction of Weyl group representations,
see Theorem 5.4. (Since $\cu_G^{sp}$ is naturally in bijection with the set of
special representations of the Weyl group $W$,
we formulate our results in terms of $W$. This has the advantage
that our results make sense for any finite Coxeter group.)

\subhead 0.2\endsubhead
In the remainder of this paper we fix a finite Coxeter group $W$.
Let $\Irr(W)$ be the set of isomorphism classes of
irreducible representations of $W$ over $\CC$. 
In \cite{L79} a certain subset $\cs_W$ of $\Irr(W)$ was defined assuming that
$W$ is a Weyl group, but the same definition can be given in general. (Later,
the name of ``special representations'' was given to the elements of this
subset.) Writing $W=\prod_iW_i$ where $W_i$ are irreducible Coxeter groups we have
that $\cs_W$ consists of all $E=\boxtimes_iE_i$ with $E_i\in\cs_{W_i}$ for all
$i$. Let $\sgn$ be the sign representation of $W$. Let
$\cs_W^{odd}=\{E\in\cs_W;E\ot\sgn\n\cs_W\}$.
When $W$ is irreducible of type $E_7/E_8/H_3/H_4$, $\cs_W^{odd}$ consists of
the $E\in\cs_W$ which have dimension $512/4096/4/16$ (there are $1/2/1/2$ such $E$); when $W$ is irreducible of type other than $E_7,E_8,H_3,H_4$, then
$\cs_W^{odd}=\emp$.
We set $\cs_W^{ord}=\cs_W-\cs_W^{odd}$.
If $E\in\cs_W$ we say that $E$ is odd if $E\in\cs_W^{odd}$; 
we say that $E$ is ordinary if $E\in\cs_W^{ord}$.

When
$W$ is irreducible let $E\m E^\cir$ be the involution of $\cs_W$
such that $E^\cir=E\ot\sgn$ for $E\in\cs_W^{ord}$ and $E\m E^\cir$
is the permutation of $\cs_W^{odd}$ of order $1/2/1/2$ if $W$ is of type
$E_7/E_8/H_3/H_4$. When $W$ is not necessarily irreducible and $W=\prod_iW_i$
with $W_i$ irreducible then the involution $E\m E^\cir$ of $\cs_W$ is defined
by $\boxtimes_iE_i\m\boxtimes_i(E_i^\circ)$ where $E_i\in\cs_{W_i}$.

We can identify $\cs_W$ with the set $Cell_W$
of two-sided cells of $W$ by the procedure stated in \cite{KL, 1.7}. 
Then the involution $E\m E^\cir$ of $\cs_W$ becomes the
involution of $Cell(W)$ given by left multiplication
by the longest element of $W$; the partial order $\le_{LR}$ on $Cell_W$
in \cite{KL} becomes a partial order on $\cs_W$ denoted by $\le$.
Note that $\sgn$ is the unique maximal element for $\le$. Using
\cite{KL, 3.3} we see that

(a) $E\le E'$ (in $\cs_W$) implies $E'{}^\cir\le E^\cir$.
\nl
Let $Cell_W@>>>\NN$ be the function defined in \cite{L85, 5.4} (the definition
given in {\it loc.cit.} for $W$ assumed to be a Weyl group is also applicable
without this assumption). This can be viewed as a function $E\m a_E,\cs_W@>>>\NN$.
Using {\it loc.cit.} we see that

(b) $E\le E'$ (in $\cs_W$) implies $a_{E'}\le a_E$.
\nl
Let $\hcs_W$ be the set of all $(E,E')\in\cs_W\T\cs_W$
such that $E\lneqq E'$ and there is no $E''\in\cs_W$ with
$E\lneqq E'',E''\lneqq E'$.
For $E,E'$ in $\cs_W$ we say that $E,E'$ are
{\it adjacent} if $(E,E')\in\hcs_W$ or $(E',E)\in\hcs_W$;
we then write $E---E'$.
Note that

(c) $E---E'$ implies $E^\cir---E'{}^\cir$.
\nl
(See (a).)

Writing $W=\prod_iW_i$ where $W_i$ are irreducible Coxeter groups we have
that $E---E'$ for $E=\boxtimes_iE_i,E'=\boxtimes_iE'_i$ with $E_i,E'_i$ in $\cs_{W_i}$ for all
$i$ if and only if $E_i---E'_i$ (in $\cs_{W_i}$) for one $i$ and $E_i=E'_i$ for all other $i$.

\proclaim{Theorem 0.3}Assume that $W$ is irreducible.
Let $(E,E')$ be in $\hcs_W$.

(a) We have either $a_E-a_{E'}=1$ or $a_{E'{}^\cir}-a_{E^\cir}=1$.

(b) To $E,E'$ one can associate an irreducible Coxeter group $W'=W_{E,E'}$
contained in $W$ such that $a_E-a_{E'}=h-1$ where $h$ is the Coxeter number of $W_{E,E'}$.
\endproclaim

\subhead 0.4\endsubhead
Let $E,E',h$ be as in 0.3(b). The condition in 0.3(b) does not define $W'=W_{E,E'}$ uniquely
(unless $h=2$ or $h=3$). Nevertheless I believe that in each case there is a natural choice for
$W'$. We describe it below.

If $W$ is of type $A_n,n\ge1$ then $W'$ is of type $A_{h-1}$.
If $W$ is of type $B_n,n\ge2$ then $h$ is even and $W'$ is of type $B_{h/2}$.
If $W$ is of type $D_n,n\ge4$ then $h$ is even and $W'$ is of type $B_{h/2}$ or
of type $D_{(h+2)/2}$; more precisely, if $E$ (resp. $E'$) is associated to the
pair of partitions
$1^{n-j},1^j$ (resp. $1^{n-j-1},1^{j+1}$) of $n$ with $0\le j<(n-2)/2$,
so that $h=2n-2-4j$, then $W'$ is of type $D_{(h+2)/2}=D_{n-2j}$; in all other cases
$W'$ is of type $B_{h/2}$.

If $W$ is a dihedral group of order $2n\ge8$, then $h\in\{2,n\}$;
if $h=2$ then $W'$ is of type $A_1$; if $h=n$ then $W'=W$.

If $W$ is of type $E_8,E_7,E_6,F_4,G_2,H_3,H_4$, the various $E_1,E'_1$ such that \lb
$E_1---E'_1$ are listed in \S4 in the form
$E_1\overset{W'}\to{---}E'_1$ where $W'=W_{E_1,E'_1}$; in the case where
$a_{E_1}-a_{E'_1}=\pm1$ we omit writing $W'$ (it is of type $A_1$).
In these cases we denote elements $E_1\in\cs_W$ in the form
$d_a$ where $d=\dim E_1,a=a_{E_1}$; when $W$ is of type $F_4$
there may be two $E_1$ with the same $d,a$; we denote them by
$d_a,d'_a$.

\head 1. The polynomials $\ti\Pi_{E,E'}$\endhead
\subhead 1.1\endsubhead
We now assume (until the end of 1.4) that $W$ is the Weyl group of a reductive connected
group $G$ over $\CC$. Let $\cu_G$ be as in 0.1 and let
$\tcu_G$ be the set of pairs $(C,\cl)$ where $C\in\cu_G$ and $\cl$ is an
irreducible local system on $C$, equivariant for conjugation by $G$. Springer's
correspondence provides an imbedding $\io:\Irr(W)@>>>\tcu_G$.
As stated in \cite{L79, no.9}, there is a well defined subset $\cu_G^{sp}$ of $\cu_G$
such that $\io$ defines a bijection of $\cs_W$ onto
$\{(C,\cl)\in\tcu_G,C\in\cu_G^{sp},\cl=\CC\}$.
This can be viewed as a bijection $\io':\cs_W@>\si>>\cu_G^{sp}$.
Let $\le$ be the partial order on $\cu_G^{sp}$ obtained by restricting the obvious
partial order on $\cu_G$ (the unique maximal element is the regular unipotent class).
According to \cite{B09} or \cite{G12}, $\io'$ is compatible with the partial
orders on $\cs_W,\cu_G^{sp}$. It follows that the order reversing involution
$E\m E^\cir$ on $\cs_W$ can be viewed as an order reversing involution of $\cu_G^{sp}$.
Now the partial order on $\cu_G$ has been explicitly computed in all cases.
(See \cite{Sp} and the references there.)
Since $\io,\io'$ are explicitly known, this determines the partial order (hence the
adjacency relation) on $\cs_W$. In the case where $W$ is of type $E_8,E_7,E_6,F_4,G_2$
the adjacency relation on $\cs_W$ is described in \S4;
from this one can see that 0.3(a) holds in these cases.

\subhead 1.2\endsubhead
Let $u$ be an indeterminate.
Let $E,E'$ be in $\Irr(W)$. Let
$(C,\cl)=\io(E)$, $(C',\cl')=\io(E')$. Assume that $C\sub \bC'$ where $\bC'$ is the
closure of $C'$. For $i\in\ZZ$ let
$n_{E,E',i}$ be the multiplicity of $\cl$ in the restriction to $C$ of the
$i$-th cohomology sheaf of the intersection cohomology complex of $\bC'$ with coefficients in $\cl'$.
It is known that $n_{E,E',i}$ is zero unless $i\in2\NN$; we set
$\Pi_{E,E'}=\sum_{i\in2\NN}n_{E,E',i}u^{i/2}\in\NN[u]$.

\subhead 1.3\endsubhead
Assume that $W$ is irreducible. Assume that $(E,E')$ is in $\hcs_W$.
If $E$ and $E'$ are ordinary (see 0.2) we set $\ti\Pi_{E,E'}=\Pi_{E,E'}$.
If $E$ is ordinary and $E'$ is odd (see 0.2) we set
$\ti\Pi_{E,E'}=\Pi_{E,E'}+\Pi_{E,E'_1}$
where $E'_1\in\Irr(W)$ is in the same two-sided cell as $E'$ and is distinct from $E'$.
If $E$ is odd and $E'$ is ordinary we set
$\ti\Pi_{E,E'}=\Pi_{E,E'}+\Pi_{E_1,E'}$
where $E_1\in\Irr(W)$ is in in the same two-sided cell as $E$ and is distinct from $E$.
(Note that at most one of $E,E'$ is odd). 

\proclaim{Conjecture 1.4} In the setup of 1.3 we have
$$\ti\Pi_{E,E'}=u^{e_1-1}+u^{e_2-1}+\do+u^{e_r-1}$$
where $e_1\le e_2\le\do\le e_r$ are the exponents of a well defined
irreducible Weyl group $W'$ of rank $r\ge1$.
\endproclaim
Using the known tables (in particular those of \cite{BS} for type $E_6,E_7,E_8$)
I have verified this for $W$ of type
$B_2,B_3,D_4,D_5,D_6,D_7,E_6,E_7,E_8,F_4,G_2$. The resulting $W'$ is
described for types $E_8,E_7,E_6,F_4,G_2$ in \S4.
We can thus verify 0.3(b) (in the strengthened form 0.4) in these cases.
When $W$ is of type $A$ the conjecture can be deduced from \cite{KP81}; this
also verifies 0.3(a) and 0.3(b) (in the strengthened form 0.4) in this case; an alternative
proof is given in \S3.

\subhead 1.5\endsubhead
We now assume that $W$ is irreducible but not a Weyl group.
In this case the partial order $\le$ on $\cs_W$ is linear (it is
described in \cite{G12}); the adjacency relation is easily described
and 0.3(a) is easily verified (see \S4 for types $H_4,H_3$).
For any $E,E'$ in $\Irr(W)$ we set
$\Pi_{E,E'}=u^{-a_{E'}}P_{E',E}$ where $P_{E',E}\in\QQ(u)$ is
attached by Geck and Malle \cite{GM} to $E,E'$.
It is known that in our case we have $\Pi_{E,E'}\in\NN[u]$.
Assuming that $(E,E')\in\hcs_W$,
we modify $\Pi_{E,E'}$ to a $\ti\Pi_{E,E'}$
by the same procedure as in 1.3.
Then conjecture 1.4 can be extended word by word to our case
(except that $W'$ is now a finite irreducible Coxeter group, not
necessarily a Weyl group). This conjecture is actually true, as can
be verified using the tables for $P_{E',E}$ available through CHEVIE.
The values of $W'$ for $W$ of type $H_4,H_3$ are listed in \S4.
We can thus verify 0.3(b) (in the strengthened form 0.4) in these cases.
In the case where $W$ is a finite dihedral group. 0.3(a) and 0.3(b)
(in the strengthened form 0.4) are easily verified.

If $W$ is a Weyl group, it is likely that the procedure above (based on \cite{GM})
gives the same polynomials as in 1.3.

\head 2. Type B,D\endhead
\subhead 2.1\endsubhead
Let $\e\in\{0,1\}$. Let $r\ge1$.  For $M\in\e+2\NN$ let
$\fS_r^M$ be the set of all $(a_1,a_2,...,a_M)\in\NN^M$ such that $a_1\le a_2\le\do\le a_M$
with no two consecutive equal signs and such that $a_1+a_2+...+a_M=r+(M^2-2M+\e)/4$.
For example when $\e=1$ we have $(0,1,1,2,2,...,r,r)\in\fS_r^{2r+1}$, $(r)\in\fS_r^1$;
when $\e=0$ we have $(0,1,1,2,2,...,r-1,r-1,r)\in\fS_r^{2r}$, $(0,r)\in\fS_r^2$.

We define an equivalence relation on $\fS_r^\e\cup\fS_r^{\e+2}\cup\fS_r^{\e+4}\cup\do$ by
$$\align&(a_1,a_2,\do,a_M)\si(0,0,a_1+1,a_2+1,\do,a_M+1)\\&
\si(0,0,1,1,a_1+2,a_2+2,\do,a_M+2)\si\do.\endalign$$
Let ${}^\e\fS_r$ be the set of equivalence classes.

\subhead 2.2\endsubhead
If $a_*=(a_1,a_2,...,a_M)\in\fS_r^M$ and $t\ge a_M$ let $a_*^{!t}$ be
the sequence obtained from 
$0,0,1,1,2,2,...,t,t$ by removing two entries $t-a$ for any $a$ which appears twice in $a_*$
and by removing one entry $t-a$  for any $a$ which appears exactly once in $a_*$.
Now $a_*^{!t}$ has $M'=2(t+1)-M$ entries whose sum is
$$t(t+1)-tM+r+(M^2-2M+\e)/4=r+((2t+2-M)^2-2(2t+2-M)+\e)/4$$
so that $a_*^{!t}\in\fS_r^{M'}$. It is easy to see that the equivalence class of
$a_*^{!t}\in\fS_r^{2(t+1)-M}$ is independent of $t$ (if $t\ge a_M$). From the definition
we see that if $a_*$ is replaced by an equivalent sequence, the equivalence class of
$a_*^{!t}$ (with $t\ge a_M$) is unchanged. Hence $a_*\m a_*^{!t}$ defines a map $\x\m\x^\bul$,
${}^\e\fS_r\to{}^\e\fS_r$. Its square is 1.

\subhead 2.3\endsubhead
For any $p\ge2$ let $({}^\e\fS_r\T{}^\e\fS_r)_p$ be the subset of ${}^\e\fS_r\T{}^\e\fS_r$
consisting of pairs
$(\x,\x')$ which can be represented by $(a_*,a'_*)\in\fS_r^M\T\fS_r^M$ such that
for some $k<k+1$ in $[1,M]$ and some $a$ we have
$$(a_k,a_{k+1})=(a,a+p), (a'_k,a'_{k+1})=(a+1,a+p-1), a_s=a'_s \text{ for }s\ne k,k+1.$$
For any $p\ge3$ let $({}^\e\fS_r\T{}^\e\fS_r)'_p$ be the subset of ${}^\e\fS_r\T{}^\e\fS_r$
consisting of pairs
$(\x,\x')$ which can be represented by $(a_*,a'_*)\in\fS_r^M\T\fS_r^M$ such that
for some $k<l$ in $[1,M]$ with $l=k+2p-3$ and some $a$ we have
$$\align&(a_k,a_{k+1},\do,a_{l-1},a_l)\\&
=(a,a+1,a+2,a+2,a+3,a+3,\do,a+p-2,a+p-2,a+p-1,a+p),\endalign$$
$$\align&(a'_k,a'_{k+1},\do,a'_{l-1},a'_l)\\&
=(a+1,a+1,a+2,a+2,\do,a+p-2,a+p-2,a+p-1,a+p-1),\endalign$$
and $a_s=a'_s$ for $s\n\{k,k+1,\do,l-1,l\}$.

From the definitions we see that $(\x,\x')\m(\x'{}^\bul,\x^\bul)$ is a bijection
$$({}^\e\fS_r\T{}^\e\fS_r)_p@>\si>>({}^\e\fS_r\T{}^\e\fS_r)'_p$$
(if $p\ge3$) and is a bijection
$$({}^\e\fS_r\T{}^\e\fS_r)_2@>\si>>({}^\e\fS_r\T{}^\e\fS_r)_2.$$

\subhead 2.4\endsubhead
When $W$ is of type $B_r$, $r\ge1$, we
can identify $\cs_W$ with ${}^1\fS_r$ as in \cite{L79, no.5}.
When $W$ is of type $D_r$, $r\ge2$, we
can identify $\cs_W$ with ${}^0\fS_r$ as in \cite{L79, no.5}
except that each element of ${}^0\fS_r$ of the form
$(a_1,a_2,...,a_M)$ with $a_1=a_2<a_3=a_4<a_5=a_6<\do$
corresponds to two representations in $\cs_W$.
(We use the notation $(a_1,a_2,a_3,a_4,\do)$ instead of the notation
$\left(\smallmatrix a_1&a_3&a_5&\do\\a_2&a_4&a_6&\do\esm\right)$ in \cite{L79}.)
If $E,E'$ in $\cs_W$ correspond to $\x,\x'$ in ${}^\e\fS_r$ then $E,E'$ are adjacent in $\cs_W$
if and only if (i) or (ii) below holds:

(i) $(\x,\x')$ or $(\x',\x)$ belongs to $({}^\e\fS_r\T{}^\e\fS_r)_p$ for some $p\ge2$;

(ii) $(\x,\x')$ or $(\x',\x)$ belongs to $({}^\e\fS_r\T{}^\e\fS_r)'_p$ for some $p\ge3$.
\nl
If (i) holds then $a_E-a_{E'}=\pm1$; if (ii) holds then $a_E-a_{E'}=\pm(2p-3)$. The involution
$E\m E^\cir$ of $\cs_W$ corresponds to the involution $\x\m\x^\bul$ of
${}^\e\fS_r$ and that involution interchanges pairs as in (i) with pairs as in (ii) (if $p\ge3$)
or pairs in (i) with pairs in (i) (if $p=2$); see 2.3. It follows that 0.3 holds in our case.

\head 3. Type $A$\endhead
\subhead 3.1\endsubhead
Let $r\ge1$. For $m\in\NN$ let
$\fT_r^m$ be the set of all $(a_1,a_2,\do,a_m)\in\NN^m$ such that $a_1<a_2<\do<a_m$
and such that $a_1+a_2+\do+a_m=r+m(m-1)/2$.
For example we have $(1,2,\do,r)\in\fT_r^r$, $(r)\in\fT_r^1$.
We define an equivalence relation on $\fT_r^0\cup\fT_r^1\cup\fT_r^2\cup\do$ by
$$\align&(a_1,a_2,\do,a_m)\si(0,a_1+1,a_2+1,\do,a_m+1)\\&
\si(0,1,a_1+2,a_2+2,\do,a_m+2)\si\do.\endalign$$
Let $\fT_r$ be the set of equivalence classes.

\subhead 3.2\endsubhead
If $a_*=(a_1,a_2,\do,a_m)\in\fT_r^m$ and $t\ge a_m$ let $a_*^{!t}$ be the sequence obtained
from  $0,1,2,\do,t$ by removing an entry $t-a$ for any $a$ which appears in $a_*$.
Now $a_*^{!t}$ has $m'=t+1-m$ entries whose sum is
$$t(t+1)/2-tm+r+m(m-1)/2=r+(t+1-m)(t-m)/2$$
so that $a_*^{!t}\in\fT_r^{m'}$. It is easy to see that the equivalence class of
$a_*^{!t}\in\fT_r^{t+1-m}$ is independent of $t$ (if $t\ge a_m$). 
From the definition we see that if $a_*$ is replaced by an equivalent sequence,
the equivalence class of $a_*^{!t}$ (with $t\ge a_m$) is unchanged. Hence $a_*\m a_*^{!t}$ defines a map $\x\m\x^\bul$,
$\fT_r\to\fT_r$. Its square is 1.

\subhead 3.3\endsubhead
For any $p\ge3$ let $(\fT_r\T\fT_r)_p$ be the subset of $\fT_r\T\fT_r$
consisting of pairs $(\x,\x')$ which can be represented by $(a_*,a'_*)\in\fT_r^m\T\fT_r^m$
such that
for some $k<k+1$ in $[1,m]$ and some $a$ we have
$$(a_k,a_{k+1})=(a,a+p), (a'_k,a'_{k+1})=(a+1,a+p-1), a_s=a'_s \text{ for }s\ne k,k+1.$$
For any $p\ge4$ let $(\fT_r\T\fT_r)'_p$ be the subset of $\fT_r\T\fT_r$
consisting of pairs $(\x,\x')$ which can be represented by $(a_*,a'_*)\in\fT_r^m\T\fT_r^m$
such that for some $k<l$ in $[1,m]$ with $l=k+p-1$ and some $a$ we have
$$(a_k,a_{k+1},\do,a_{l-1},a_l)    =(a,a+2,a+3,\do,a+p-3,a+p-2,a+p),$$
$$(a'_k,a'_{k+1},\do,a'_{l-1},a'_l)=(a+1,a+2,a+3,\do,a+p-3,a+p-2,a+p-1),$$
and $a_s=a'_s$ for $s\n\{k,k+1,...,l-1,l\}$.

From the definitions we see that $(\x,\x')\m(\x'{}^\bul,\x^\bul)$ is a bijection
$$(\fT_r\T\fT_r)_p@>\si>>(\fT_r\T\fT_r)'_p$$
(if $p\ge4$) and is a bijection
$$(\fT_r\T\fT_r)_3@>\si>>(\fT_r\T\fT_r)_3.$$

\subhead 3.4\endsubhead
When $W$ is of type $A_{r-1}$, $r\ge2$,
we can identify $\cs_W$ with $\fT_r$ in the standard way;
the unit representation corresponds to $(1,2,\do,r)$; $\sgn$ corresponds to $(r)$.

If $E,E'$ in $\cs_W$ correspond to $\x,\x'$ in $\fT_r$ then $E,E'$ are adjacent in $\cs_W$
if and only if (i) or (ii) below holds:

(i) $(\x,\x')$ or $(\x',\x)$ belongs to $(\fT_r\T\fT_r)_p$ for some $p\ge3$;

(ii) $(\x,\x')$ or $(\x',\x)$ belongs to $(\fT_r\T\fT_r)'_p$ for some $p\ge4$.
\nl
If (i) holds then $a_E-a_{E'}=\pm1$; if (ii) holds then $a_E-a_{E'}=\pm(p-2)$.
The involution $E\m E^\cir$ of $\cs_W$ corresponds to the involution $\x\m\x^\bul$ of
$\fT_r$ and that involution interchanges pairs as in (i) with pairs as in (ii) (if $p\ge4$)
or pairs in (i) with pairs in (i) (if $p=3$). It follows that 0.3 holds in our case.

\head 4.Examples\endhead
\subhead 4.1\endsubhead
In this section we assume that $W$ is of type $E_8,E_7,E_6,F_4,G_2,H_3,H_4$. In these
cases we describe the pairs $(E,E')$ of adjacent elements of $\cs_W$
in the form $E\overset{W'}\to{---}E'$ where $W'=W_{E,E'}$ is determined by 1.4, 1.5 (when
$a_E-a_{E'}=\pm1$ we omit writing $W'$).

\subhead 4.2. Type $E_8$\endsubhead

$35_2---8_1---1_0$,  $1_{120}\overset{E_8}\to{---}
8_{91}\overset{E_7}\to{---}35_{74}$

$210_4---112_3---35_2$,  $35_{74}\overset{B_6}\to{---}112_{63}\overset{F_4}\to{---}210_{52}$

$567_6---560_5---210_4$,  $210_{52}\overset{B_3}\to{---}560_{47}---567_{46}$

$1400_7---700_6---560_5$,  $560_{47}\overset{B_3}\to{---}700_{42}\overset{G_2}\to{---}1400_{37}$

$525_{12}\overset{G_2}\to{---}1400_7---567_5$,  $567_{46}\overset{G_2}\to{---}1400_{37}---525_{36}$

$2240_{10}---3240_9---1400_8---1400_7$,

$1400_{37}\overset{G_2}\to{---}1400_{32}---3240_{31}\overset{B_2}\to{---}2240_{18}$

$4096_{11}---2268_{10}---3240_9$, $3240_{31}---2268_{30}\overset{A_4}\to{---}4096_{26}$

$2800_{13}\overset{A_2}\to{---}4096_{11}---2240_{10}$,
$2240_{28}\overset{A_2}\to{---}4096_{26}---2800_{25}$

$6075_{14}---2800_{13}---525_{12}$,  $525_{36}\overset{F_4}\to{---}2800_{25}\overset{B_2}\to{---}6075_{22}$

$4536_{13}---4200_{12}---4096_{11}$, $4096_{26}\overset{A_2}\to{---}4200_{24}---4536_{23}$

$4480_{16}---5600_{15}---6075_{14}---4536_{13}$,

$4536_{23}---6075_{22}---5600_{21}\overset{G_2}\to{---}4480_{16}$

$4200_{15}---2835_{14}---4536_{13}$,
$4536_{23}---2835_{22}---4200_{21}$

$4480_{16}---4200_{15}---6075_{14}$,  $6075_{22}---4200_{21}\overset{G_2}\to{---}4480_{16}$

$2100_{20}\overset{B_3}\to{---}5600_{15}$, $5600_{21}---2100_{20}$

\subhead 4.3. Type $E_7$\endsubhead

$27_2---7_1---1_0$,  $1_{63}\overset{E_7}\to{---}7_{46}\overset{D_6}\to{---}27_{37}$

$120_4---56_3---27_2$,  $27_{37}\overset{B_4}\to{---}56_{30}\overset{B_3}\to{---} 120_{25}$

$120_4---21_3---27_2$,  $27_{37}---21_{36}\overset{F_4}
\to{---} 120_{25}$

$105_6---189_5---120_4$,  $120_{25}\overset{B_2}\to{---}189_{22}---105_{21}$

$315_7---168_6--- 189_5$,  $189_{22}---168_{21}\overset{G_2}\to{---}315_{16}$

$315_7---210_6--- 189_5$,  $189_{22}---210_{21}\overset{B_3}\to{---}315_{16}$

$105_{12}---315_7---105_6$, $105_{21}\overset{G_2}\to{---}315_{16}---105_{15}$

$168_6---189_7---210_6$, $210_{21}---189_{20}---168_{21}$

$378_9---405_8---189_7$,
$189_{20}\overset{B_3}\to{---}405_{15}---378_{14}$

$105_{15}\overset{B_3}\to{---}420_{10}---378_9$,
$378_{14}---420_{13}---105_{12}$

$512_{11}---210_{10}---378_9$,
$378_{14}---210_{13}\overset{A_2}\to{---} 512_{11}$

$512_{11}---420_{10}$,
$420_{13}\overset{A_2}\to{---}512_{11}$

\subhead 4.4. Type $E_6$\endsubhead

$20_2---6_1---1_0$,  $1_{36}\overset{E_6}\to{---}6_{25}
\overset{A_5}\to{---}20_{20}$

$64_4---30_3---20_2$,  $20_{20}\overset{B_3}\to{---}
30_{15}\overset{A_2}\to{---} 64_{13}$

$80_7---81_6---60_5---64_4$, $64_{13}\overset{A_2}\to
{---}60_{11}---81_{10}\overset{B_2}\to{---} 80_7$

$80_7---24_6\overset{A_2}\to{---}64_4$,
$64_{13}---24_{12}\overset{G_2}\to{---}80_7$

\subhead 4.5. Type $F_4$\endsubhead

$9_2---4_1---1_0$,  $1_{24}\overset{F_4}\to{---}4_{13}
\overset{B_2}\to{---}9_{10}$

$12_4---8_3---9_2$,
$9_{10}---8_9\overset{G_2}\to{---}12_4$

$12_4---8'_3---9_2$,
$9_{10}---8'_9\overset{G_2}\to{---}12_4$

\subhead 4.6. Type $G_2$\endsubhead

$2_1---1_0$     $1_6\overset{G_2}\to{---}2_1$

\subhead 4.7. Type $H_4$\endsubhead

Let $H_2$ (resp. $\D$) be  a dihedral group of order $10$ (resp. $20$).

$9_2---4_1---1_0$,
$1_{60}\overset{H_4}\to{---}4_{31}\overset{\D}\to
{---}9_{22}$

$25_4---16_3---9_2$,
     $9_{22}\overset{H_2}\to{---}16_{18}
     \overset{A_2}\to{---}25_{16}$

$24_6---36_5---25_4$,
$25_{16}---36_{15}\overset{\D}\to{---}24_6$

\subhead 4.8. Type $H_3$\endsubhead

$5_2---3_1---1_0$,
$1_{15}\overset{H_3}\to{---}3_6---5_5$

$ 4_3---5_2$,   $5_5\overset{A_2}\to{---}4_3$

\head 5. Truncated induction and adjacency\endhead
\subhead 5.1\endsubhead
In this section we assume that $W$ is a Weyl group.
Let $W_J$ be a standard parabolic subgroup of $W$. In \cite{LS} a function
$j_{W_J}^W$ (``truncated induction'')
from a certain subset of $\Irr(W_J)$ to $\Irr(W)$ was defined.
From \cite{L79} we see that $\cs_{W_J}$ is contained in the subset on which
$j_{W_J}^W$ is defined and that the image of
$\cs_{W_J}$ under $j_{W_J}^W$ is contained in $\cs_W$. Thus we have a well defined
map $j_{W_J}^W:\cs_{W_J}@>>>\cs_W$. This map preserves the value of the $a$-function.

\subhead 5.2\endsubhead
The following result appears in \cite{L79, no.6}. 

(a) {\it Assume that $W$ is irreducible, $\ne1$. Let $E\in\cs_W$. Then either $E$ or $E^\cir$ is
of the form $j_{W_J}^W(E_1)$ where $W_J\sub W$ is as in 5.1 with $W_J\ne W$ and
$E_1\in\cs_{W_J}$.}

Thus $\cs_W$ admits an inductive definition based on the involution $E\m E^\cir$ and on
truncated induction from proper parabolic subgroups.
In \cite{G12} a definition of $\le$ on $\cs_W$ in the same spirit (but with the full
induction playing a role) was given.

\subhead 5.3\endsubhead
We want to show that the adjacency relation  on $\cs_W$ can be described in terms of
$E\m E^\cir$ and truncated induction.

We define a subset $\check{\cs}_W$ by induction on the rank of $W$.
Writing $W=\prod_iW_i$ where $W_i$ are irreducible Weyl groups we have
that $(E,E')\in\check{\cs}_W$ for $E=\boxtimes_iE_i,E'=\boxtimes_iE'_i$ with $E_i,E'_i$ in
$\cs_{W_i}$ for all $i$ if and only if $(E_i,E'_i)\in\check{\cs}_{W_i}$ for one $i$ and
$E_i=E'_i$ for all other $i$.
If $W=\{1\}$, we have $\check{\cs}_W=\emp$. If $W$ is of type $A_1$ then 
$\check{\cs}_W$ consists of $(\sgn,1)$.
Assume now that $W$ is irreducible of rank $\ge2$.

We say that $(E,E')\in\cs_W\T\cs_W$ is {\it induced } if
there exists $W_J\sub W$ as in 5.1 with $W_J\ne W$ and $(E_1,E'_1)\in\check{\cs}_{W_J}$
such that $a_{E_1}-a_{E'_1}=1$ ($a$-function relative to $W_J$) and $E=j_{W_J}^W(E_1)$,
$E'=j_{W_J}^W(E'_1)$.

Then $\check{\cs}_W$ consists of all
$(E,E')\in\cs_W\T\cs_W$ such that at least one of (i),(ii),(iii),(iv) below are satisfied.

(i) $(E,E')$ is induced;

(ii) $(E'{}^\cir,E^\cir)$ is induced;

(iii) $W$ is of type $B_3$ and $(E,E')$ satisfies $a_E=3,a_E=2$ (so that
$(E,E')=(E'{}^\cir,E^\cir)$);

(iv) $W$ is of type $E_7$ and $(E,E')=(512_{11},210_{10})$ (notation of 4.3).
\nl
This completes the inductive definition of $\check{\cs}_W$.

We state the following result.

\proclaim{Theorem 5.4} We have $\check{\cs}_W=\hcs_W$.
\endproclaim

\subhead 5.5\endsubhead
We show that

(a) $\check{\cs}_W\sub\hcs_W$
\nl
by induction on the rank of $W$. We can assume that
$W$ is irreducible of rank $\ge2$. Let $(E,E')\in\check{\cs}_W$.

Assume first that $(E,E')$ is induced (from $(E_1,E'_1)\in\check{\cs}_{W_J}$).
Let $C,C'$ in $\cu_G^{sp}$ correspond to $E,E'$ as in 1.1; similarly
let $C_1,C'_1$ in $\cu_{G_J}^{sp}$ correspond to $E_1,E'_1$ where $G_J$ is a Levi subgroup
of a parabolic subgroup of $G$.
Since $a_{E_1}-a_{E'_1}=1$ we see that $\dim(C'_1)-\dim(C_1)=2$; moreover $C_1$ is contained
in the closure of $C'_1$. By \cite{LS}, $C$ and $C'$ are induced from $C_1$ and $C'_1$, so that
$C$ is contained in the closure of $C'$ and $\dim(C')-\dim(C)=2$; it follows that
$E\lneqq E'$ and there is no $E''\in\cs_W$ such that $E\lneqq E'',E''\lneqq E'$
and
$(E,E')\in\hcs_W$.

Next we assume that $(E'{}^\cir,E^\cir)$ is induced. By the earlier argument we see that
$(E'{}^\cir,E^\cir)\in\hcs_W$. This implies that $(E,E')\in\hcs_W$, by 0.2(c).
If $(E,E')$ is as in 5.3(iv) then $(E,E')\in\hcs_W$ follows from 4.3.
If $(E,E')$ is as in 5.3(iii) then $(E,E')\in\hcs_W$ can be verified directly.
This proves (a).

\subhead 5.6\endsubhead
Let $\e\in\{0,1\}$. In this subsection we verify some properties of a pair
$(a_*,a'_*)\in\fS^M_r\T\fS^M_r$ with $M\in\e+2\NN$ (notation in \S2)
which will be used in the proof of the converse to 5.5(a) in type $B_r,D_r$.

We say that $(a_*,a'_*)$ is induced if there exists
$k\in\{1,2,...,r\}$ such that

$a_*=(a_1,a_2,\do,a_M)$, $a'_*=(a'_1,a'_2,\do,a'_M)$,

$(a_1,a_2,\do,a_{M-k},a_{M-k+1}-1,a_{M-k+2}-1,\do,a_M-1)\in\fS_{r-k}^M$,

$(a'_1,a'_2,\do,a'_{M-k},a'_{M-k+1}-1,a'_{M-k+2}-1,\do,a'_M-1)\in\fS_{r-k}^M$.

If $(a_*,a'_*)$ represents an element of $({}^\e\fS_r\T{}^\e\fS_r)_p$ with $p\ge3$ then
$(a_*,a'_*)$ is clearly induced. 
If $(a_*,a'_*)$ represents an element of $({}^\e\fS_r\T{}^\e\fS_r)'_p$ with $p\ge3$ and $t\gg0$
then, by 2.3, $(a'_*{}^{!t},a_*{}^{!t})$
represents an element of $({}^\e\fS_r\T{}^\e\fS_r)_p$ with $p\ge3$ hence is induced. We show:

(a) {\it Assume that $(a_*,a'_*)$ represent an element of
$({}^\e\fS_r\T{}^\e\fS_r)_p$ with $p=2$. Then}

(i) {\it $(a_*,a'_*)$ is induced or}

(ii) {\it $(a'_*{}^{!t},a_*{}^{!t})$ (with $t\gg0$) is induced or}

(iii) {\it $(a_*,a'_*)$ is equivalent to one of
$(0,2),(1,1))$ or $((0,0,2),(0,0,2))$ \lb or $((0,2,2),(1,1,2))$
or $((0,0,2,2),(0,1,1,2))$.}
\nl
We have $a_*=(\do,a,a+2,\do),a'_*=(\do,a+1,a+1,\do)$ with the same entries in the same position
marked by $\do$. We assume that $(a_*,a'_*)$ is not as in (i),(ii). 

If $a_*=(\do,a,a+2,c_1,\do),a'_*=(\do,a+1,a+1,c_1,\do)$ then $c_1=a+2$ (if $c_1\ge a+3$ then
$(a_*,a'_*)$ would be induced).

If $a_*=(\do,a,a+2,c_1,c_2,\do),a'_*=(\do,a+1,a+1,c_1,c_2,\do)$ then $c_2=a+3$
(if $c_2\ge a+4$ then $(a_*,a'_*)$ would be induced).

If $a_*=(\do,a,a+2,c_1,c_2,c_3,\do),a'_*=(\do,a+1,a+1,c_1,c_2,c_3,\do)$ then $c_3=a+3$
(if $c_3\ge a+4$ then $(a_*,a'_*)$ would be induced).

If $a_*=(\do,a,a+2,c_1,c_2,c_3,c_4,\do),a'_*=(\do,a+1,a+1,c_1,c_2,c_3,c_4,\do)$ then $c_4=a+4$
(if $c_4\ge a+5$ then $(a_*,a'_*)$ would be induced).

Continuing in this way we see that

$a_*=(\do,a,a+2,***),a'_*=(\do,a+1,a+1,***)$

where $***$ stands for the sequence 

$(a+2,a+3,a+3,a+4,a+4,a+5,a+5,\do)$

with $k\ge0$ terms.

If both $a_*,a'_*$ end with $a+s,a+s$ ($s\ge3$) then $(a'_*{}^{!(a+s)},a_*^{!(a+s)})$ is of the
the form $((b,\do),(b',\do))$ with $b>0,b'>0$ hence is induced (contradicting our assumption).

If both $a_*,a'_*$ end with $a+s-1,a+s-1,a+s$ ($s\ge4$) then
$(a'_*{}^{!(a+s)},a_*^{!(a+s)})$ is of the
form $((0,b,\do),(0,b',\do))$ with $b\ge2,b'\ge2$ hence is induced
(contradicting our assumption).

If both $a_*,a'_*$ end with $a+2,a+3$ then
$(a'_*{}^{!(a+3)},a_*^{!(a+3)})$ is of the
form $((0,b,c,\do),(0,b',c',\do))$ with $c>b\ge1,b'\ge2$, hence is induced
(contradicting our assumption).

Thus we have

(b) $(a_*,a'_*)=(\do,a,a+2),(\do,a+1,a+1)$ or

(c) $(a_*,a'_*)=(\do,a,a+2,a+2),(\do,a+1,a+1,a+2)$.
\nl
If (b) holds then $(a'_*{}^{!(a+2)},a_*^{!(a+2)})$ is of the form $((0,0,2,\do),(0,1,1,\do))$.
The first part of the argument applies to $((0,0,2,\do),(0,1,1,\do))$ and shows that it is of
the form $((0,0,2),(0,1,1))$ or of the form $((0,0,2,2),(0,1,1,2))$.

If (c) holds then 
$(a'_*{}^{!(a+2)},a_*^{!(a+2)})$ is of the form $((0,2,\do),(1,1,\do))$.
The first part of the argument applies to $((0,2,\do),(1,1,\do))$ and shows that it is of
the form $((0,2),(1,1))$ or of the form $((0,2,2),(1,1,2))$.
It follows that $(a_*,a'_*)$ is as in (iii). This proves (a).

\subhead 5.7\endsubhead
In this subsection we verify some properties of a pair
$(a_*,a'_*)\in\fT^m_r\T\fT^m_r$ with $m\in\NN$ (notation in \S3)
which will be used in the proof of the converse to 5.5(a) in type $A_{r-1}$.

We say that $(a_*,a'_*)$ is induced if there exists
$k\in\{1,2,\do,r\}$ such that

$a_*=(a_1,a_2,\do,a_m)$, $a'_*=(a'_1,a'_2,\do,a'_m)$,

$(a_1,a_2,\do,a_{m-k},a_{m-k+1}-1,a_{m-k+2}-1,\do,a_m-1)\in\fT_{r-k}^m$,

$(a'_1,a'_2,\do,a'_{m-k},a'_{m-k+1}-1,a'_{m-k+2}-1,\do,a'_m-1)\in\fT_{r-k}^m$.

If $(a_*,a'_*)$ represents an element of $(\fT_r\T\fT_r)_p$ with $p\ge4$ then
$(a_*,a'_*)$ is clearly induced. 
If $(a_*,a'_*)$ represents an element of $(\fT_r\T\fT_r)'_p$ with $p\ge4$ and $t\gg0$
then, by 3.3, $(a'_*{}^{!t},a_*{}^{!t})$ (with $t\gg0$)
represents an element of $(\fT_r\T\fT_r)_p$ with $p\ge4$ hence is induced. We show:

(a) {\it Assume that $(a_*,a'_*)$ represent an element of
$(\fT_r\T\fT_r)_p$ with $p=3$. Then}

(i) {\it $(a_*,a'_*)$ is induced or}

(ii) {\it $(a'_*{}^{!t},a_*{}^{!t})$ (with $t\gg0$) is induced or}

(iii) {\it $(a_*,a'_*)$ is equivalent to $((0,3),(1,2))$.}
\nl
We have $a_*=(\do,a,a+3,\do),a'_*=(\do,a+1,a+2,\do)$ with the same entries in the same position
marked by $\do$. We assume that $(a_*,a'_*)$ is not as in (i),(ii). 

If $a_*=(\do,a,a+3,c_1,\do),a'_*=(\do,a+1,a+2,c_1,\do)$ then $c_1=a+4$ (if $c_1\ge a+5$ then
$(a_*,a'_*)$ would be induced).

If $a_*=(\do,a,a+3,c_1,c_2,\do),a'_*=(\do,a+1,a+2,c_1,c_2,\do)$ then $c_2=a+5$
(if $c_2\ge a+6$ then $(a_*,a'_*)$ would be induced).

Continuing in this way we see that

$a_*=(\do,a,a+3,***),a'_*=(\do,a+1,a+2,***)$

where $***$ stands for the sequence 
$(a+4,a+5,a+6,\do)$ with $k\ge0$ terms.
If both $a_*,a'_*$ end with $a+s$ ($s\ge4$) then $(a'_*{}^{!(a+s)},a_*^{!(a+s)})$ is of the
the form $((b,\do),(b',\do))$ with $b>0,b'>0$ hence is induced (contradicting our assumption).
Thus we have

(b) $(a_*,a'_*)=((\do,a,a+3),(\do,a+1,a+2))$.
\nl
It follows that $(a'_*{}^{!(a+3)},a_*^{!(a+3)})$ is of the form $((0,3,\do),(1,2,\do))$.
The first part of the argument applies to $((0,3,\do),(1,2,\do))$ and shows that it is of
the form $((0,3),(1,2))$.
It follows that $(a_*,a'_*)$ is as in (iii). This proves (a).

\subhead 5.8\endsubhead
We show that

(a) $\hcs_W\sub\check{\cs}_W$.
\nl
We argue by induction on the rank of $W$. We can assume that $W$ is irreducible.
If $W$ has rank $\le1$ the result is obvious. Assume now that $W$ has rank $\ge2$.
If $W$ is of type $A_n$, $n\ge2$ then (a) follows from 5.7(a).
If $W$ is of type $B_r$, $r\ge2$ or $D_r$, $r\ge4$ then (a) follows from 5.6(a).
(Note that the only pair $(a_*,a'_*)$ which can appear in 5.6(a) is
$((0,2,2),(1,1,2))$ which corresponds to $B_3$; the other pairs in 5.6(a) appear in
lower rank.) If $W$ is of exceptional type then (a) follows from the tables in \S4 and the
explicit knowledge of truncated induction. This proves (a). Theorem 5.4 is proved.


\widestnumber\key{KP81}
\Refs
\ref\key{BS}\by W.M.Beynon and N.Spaltenstein\paper The computation of Green functions of finite Chevalley
groups of type $E_n\qua (n=6,7,8)$
\jour The University of Warwick Computer Centre Report no.23\yr1982\endref
\ref\key{B09}\by R.Bezrukavnikov\paper Perverse sheaves on affine flags and nilpotent cone of the
Langlands dual group\jour Isr. J. Math\vol170\yr2009\pages185-206\endref
\ref\key{FJLS}\by B.Fu, D.Juteau, P.Levy, E.Sommers\paper Generic singularities of nilpotent orbits
\jour Adv.in Math.\vol305\yr2017\pages1-77\endref
\ref\key{G12}\by M.Geck\paper On the Kazhdan-Lusztig order on cells and families\jour Comm. Math. Helv.
\vol87\yr2012\pages 905-927\endref
\ref\key{GM}\by M.Geck and G.Malle\paper On special pieces in the unipotent variety\jour Experim. Math.
\vol8\yr1999\pages281-290\endref
\ref\key{KL}\by D.Kazhdan and G.Lusztig\paper Representations of Coxeter groups
and Hecke algebras\jour Invent. Math.\vol53\yr1979\pages165-184\endref
\ref\key {KP81}\by H.Kraft and C.Procesi\paper Minimal singularities in $GL_n$
\jour Invent. Math\vol62\yr1981\pages503-515\endref
\ref\key {KP82}\by H.Kraft and C.Procesi\paper On the geometry of conjugacy
classes in classical groups \jour Comm. Math. Helv.\vol57\yr1982\pages 539-602\endref
\ref\key{L79}\by G.Lusztig\paper A class of irreducible representations of a Weyl
group\jour Proc. Kon. Nederl. Akad.(A)\vol82\yr1979\pages323-335\endref
\ref\key{L85}\by G.Lusztig\paper Cells in affine Weyl groups\inbook Algebraic groups and related topics
\bookinfo Adv. Stud. Pure Math.6 \publ North-Holland and Kinokuniya\yr1985\pages255-287\endref
\ref\key{LS}\by G.Lusztig and N.Spaltenstein\paper Induced unipotent classes\jour
J. Lond. Math. Soc.\vol19\yr1979\pages41-52\endref
\ref\key{Sp}\by N.Spaltenstein \book Classes unipotentes et sous-groupes de Borel\bookinfo Lect. Notes
in Math.\publ Springer Verlag\vol946\yr1982\endref
\endRefs
\enddocument